# Overlapping Circles Grid Drawn with Compass and Straightedge on an Egyptian Artifact of 14th Century BC

**Amelia Carolina Sparavigna[1] and Mauro Maria Baldi[2]**

1 Department of Applied Science and Technology, Politecnico di Torino, Torino, Italy
2 Department of Control and Computer Engineering, Politecnico di Torino, Torino, Italy
Email: amelia.sparavigna@polito.it, mauro.baldi@polito.it

## Abstract

**The study of the mathematics and geometry of ancient civilizations is a task which seems to be very difficult or even impossible to fulfil, if few written documents, or none at all, had survived from the past. However, besides the direct information that we can have from written documents, we can gain some indirect evidence on mathematics and geometry also from the analysis of the decorations we find on artifacts. Here, for instance, we will show that ancient Egyptians were able of making geometric constructions using compass and straightedge, quite before the Greek Oenopides of Chois, who lived around 450 BC, had declared some of their basic principles. In fact, a wood panel covered by an overlapping circles grid pattern, found in the tomb of Kha, an architect who served three kings of 18th Dynasty (1400-1350 BC), evidences that some simple constructions with compass and straightedge were used in ancient Egypt about nine centuries before Oenopides' time.**

**Keywords: Geometry, History of Science, Geometric Patterns, Overlapping Circles Grids.**

**Subject Areas:    History of Science, History of Geometry.**

## 1. Introduction

Porphyry of Tyre, Neoplatonic philosopher of the mid third century AD, wrote in his "Life of Pythagoras" [1], that the great philosopher of Samos, who lived in the sixth century BC, learned the mathematical sciences from the Egyptians, Chaldeans and Phoenicians. According to Porphyry, old Egyptians excelled in geometry and Phoenicians in numbers and proportions. The Chaldeans were masters of "astronomical theorems, divine rites, and worship of the Gods". In the ancient world then, besides for its spectacular civilization, Egypt was also famous for its geometry. However, it is not easy to determine precisely the specific knowledge of Egyptians. In fact, today, we possess only a limited number of texts from ancient Egypt concerning mathematics and geometry. On geometry, we have the Moscow Mathematical Papyrus and the Rhind Mathematical Papyrus. From them, we find that the ancient Egyptians knew how to compute areas of several geometric shapes and the volumes of cylinders and pyramids. The Moscow Papyrus was most likely written during the 13th dynasty, based on older material probably dating to the 12th dynasty of Egypt, roughly 1850 BC [2], whereas the Rhind Papyrus dates to around 1650 BC. It was copied by scribe Ahmes from a now-lost text from the reign of a king of the 12th dynasty. In this papyrus, we can find an approximation of number $\pi$, in a problem where a circle is approximated by an octagon [3].

Besides from these documents, we can gain some indirect information on the geometry of ancient Egypt from the analysis of their of artifacts. For instance, in [3,4], we have shown how it is possible to deduce from the geometric decorations of ancient objects some information on the value of $\pi$, approximated by a rational number. For instance, from some game disks (see Figure 1) found in the tomb of Hemaka, the chancellor of a king of the First Dynasty of Egypt, about 3000 BC, we deduced that probably Hemaka knew the value of π as 66/21=22/7=3.14.

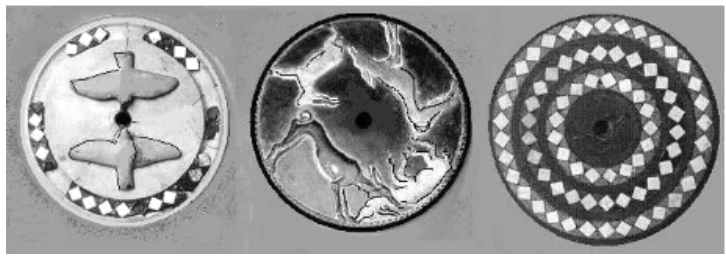

**Figure 1.** These sketches represent the gaming disks found in the tomb of chancellor Hemaka at Saqqara. Hemaka was an important official during the reign of the First Dynasty Egyptian king Den. The reign of this king lasted from 2975 BC to 2935 BC.

Here, we will show that ancient Egyptians, masters in drawing artistic figures, were also able of making constructions using compass and straightedge. It is clearly displayed by a wood panel, covered by an overlapping circles grid decoration, found in the tomb of architect Kha, who served three kings of 18th dynasty, and lived about 1400-1350 BC. Therefore, simple compass-straightedge constructions were already in use in ancient Egypt about nine centuries before the Greek Oenopides of Chois, who lived around 450 BC, had stated his principle that, for geometrical constructions, one should use no other means than compass and straightedge [5]. In fact, Oenopides' name is linked to two specific elementary constructions, one is concerning the drawing from a given point a straight line perpendicular to a given straight line, the other the drawing, on a given straight line and at a given point on it, an angle equal to a given angle [6].

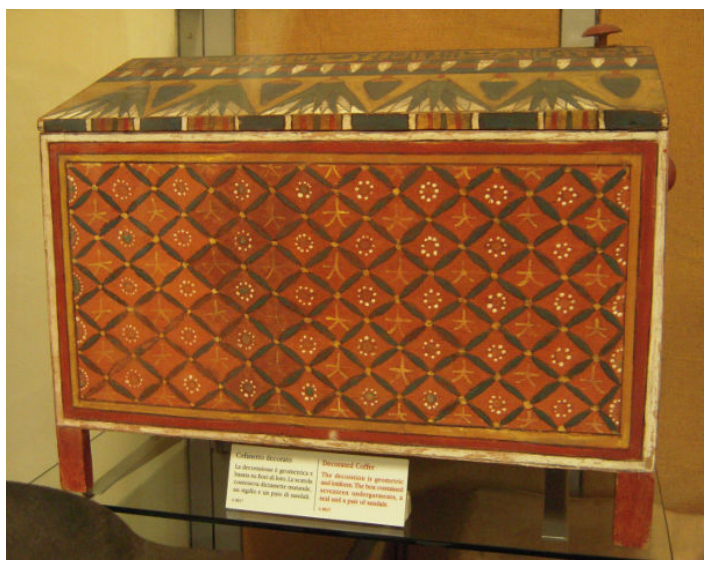

**Figure 3.** A coffin with a geometric decoration from the Kha's Tomb (Museo Egizio, Torino).

## 2. Architect Kha

Deir el Medina was a village of the ancient Egypt, inhabited by people working for the two great royal necropolises of the Valleys of the Queens and of the Kings. While excavating this village, the Italian Egyptologist Ernesto Schiaparelli discovered in February 1906 an untouched tomb [7]. It was the tomb of architect Kha, a high-ranking official of the 18th Dynasty, and his wife Merit. This tomb is rich of perfectly preserved contents and artifacts, which are showing what was the life of a wealthy family about 3400 years ago (in the Figure 2, a geometrically decorated coffin from the tomb). The contents of the tomb were transported to Turin, becoming one of the main attractions of its Museo Egizio. Among the several beautiful objects of Kha's tomb, we have discussed one in particular. It could be used as an inclinometer, a surveying instrument necessary to Kha for his work in building subterranean tombs. This device could also be the first protractor used in ancient time [8-10].

Kha was an important foreman at Deir El-Medina, responsible for building projects during the reigns of Amenhotep II, Thutmose IV and Amenhotep III [11]. These were kings of the eighteenth dynasty (c. 1543–1292 BC), the dynasty of the best known pharaohs, including Tutankhamun, Hatshepsut, the longest-reigning woman-pharaoh, and Akhenaten. The dynasty, also known as the Thutmosid Dynasty for the four kings named Thutmose, is the first of the three dynasties of the New Kingdom, the period in which ancient Egypt reached the peak of its power.

During Kha's life, Egypt was a powerful and wealthy country. The court of Amenhotep III, "was luxurious beyond imagination, with the wealth of the Mediterranean world flowing into Egypt's coffers" [12]. Therefore, Kha the architect lived in a privileged environment where, as we can easily imagine, besides the luxury and exotic goods, also knowledge and science were flowing, from the Mediterranean and Eastern countries into the cultural elite of Egypt. So, besides the skills he used for his work as architect and engineer, Kha had also received education in mathematics and geometry for sure.

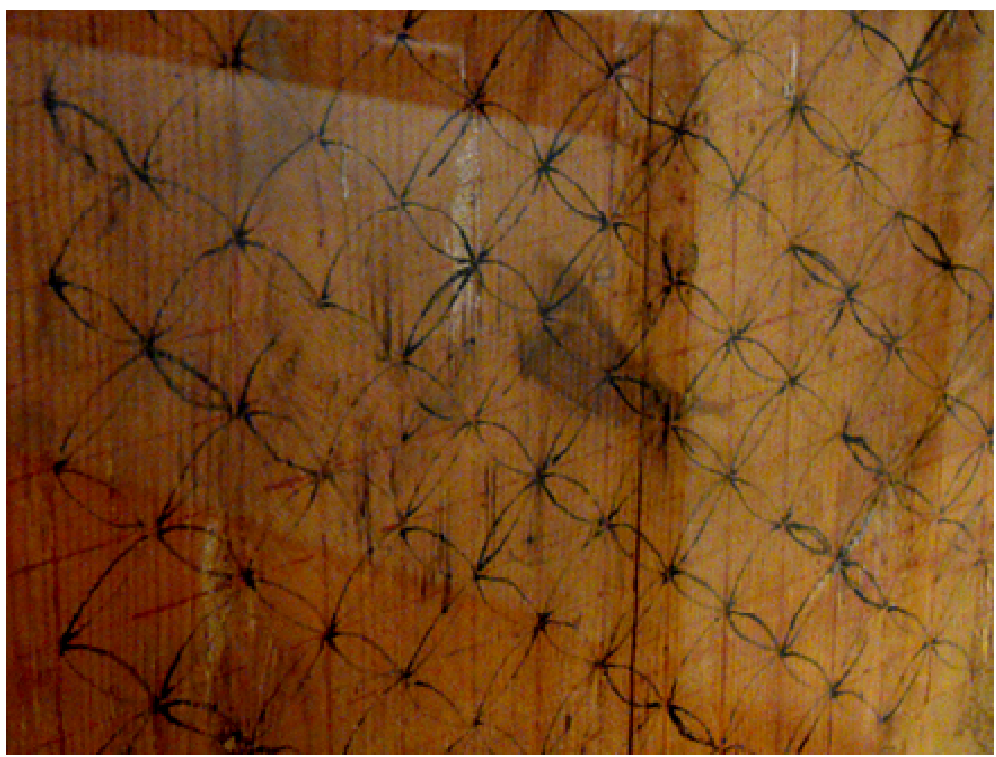

**Figure 4.** Wood panel found in the Kha's Tomb. Note the red lines creating the lattice for the centers of the black circles.

## 3. Overlapping circles grid drawn with compass and straightedge

Besides the coffins of Kha and Merit, the tomb contained all the objects necessary in their afterlife. There were the instruments that Kha used for his work and the make-up of Merit, such as ointments and kohl held in alabaster and faience vessels. Other objects in the tomb furniture include a bed, a chair, small tables and coffins, vessels and clothing garments. Geometric and polychrome figures are decorating the surface of coffins, such as that shown in Figure 2. The decoration consists of a very in-teresting and regular geometric pattern, clearly based on a regular lattice. However, it is not this artifact which is relevant to our discussion on the use in ancient Egypt of compass and straightedge. Among the objects found in the tomb, there is a wood panel which has, on its surface, a very remarkable drawing. The panel is shown in the Figure 3. Probably, it is one of the drawings which were usually made for painting on them the final polychrome decorations. It is made of parallel and perpendicular red lines and of black circles, which are giving a geometric pattern of repeating, overlapping circles of equal radii on the two-dimensional planar surface.

Several examples of polychrome decorations having such patterns were shown in the 1856 book "The Grammar of Ornament", written by Owen Jones. As told in this book, the "formation of patterns by the equal division of similar lines, as by weaving, would give to a rising people the first notions of symmetry, arrangement, disposition, and the distribution of masses". Let us add that, besides the notions of symmetry and arrangement of masses, necessary for their architecture too, the ancient Egyptians had also some notions of geometric construction with compass and straightedge. Probably, these notions were well-known when Kha was living and used for the art of decoration during the New Kingdom. A further investigation is necessary to determine if they were developed in this period or before.

## 4. Conclusion

Assuming that artists were using geometric drawings - such as that shown in the Figure 3 - for polychrome decorations such as that in the Figure 2, from the analysis of these decorations we can deduce the geometric constructions the Egyptian were able of making with compass and straightedge. From the Figure 3, for instance, it is possible to argue that they were able of drawing from a point a straight line perpendicular to a given straight line. Then, thanks to a wood panel from the tomb of architect Kha, we can tell that ancient Egyptians used, about nine centuries before, a construction with compass and straightedge described by the Greek Oenopides of Chois in the fifth century BC.